# A KIEFER–WOLFOWITZ COMPARISON THEOREM FOR WICKSELL'S PROBLEM

By Xiao Wang and Michael Woodroofe

*University of Maryland Baltimore County and University of Michigan*

We extend the isotonic analysis for Wicksell's problem to estimate a regression function, which is motivated by the problem of estimating dark matter distribution in astronomy. The main result is a version of the Kiefer–Wolfowitz theorem comparing the empirical distribution to its least concave majorant, but with a convergence rate $n^{-1} \log n$ faster than $n^{-2/3} \log n$. The main result is useful in obtaining asymptotic distributions for estimators, such as isotonic and smooth estimators.

**1. Introduction.** Our starting point is Groeneboom and Jongbloed's [5] analysis of Wicksell's [12] "Corpuscle Problem," in anatomy: Given cross-sections of a large number of corpuscles of different sizes, the distribution of radii of corpuscles was to be estimated. Assuming corpuscles were spherical, the relation between the distribution of the corpuscle radii and the distribution of the observable circular sections was derived. Wicksell approximated the density of spherical radii by a step function and then used the distributional relationship between spherical radii and circular radii to estimate the distribution of spherical radii. Groeneboom and Jongbloed [5] showed how isotonic techniques can be used in Wicksell's problem. They related the distribution of spherical radii to a nonincreasing function that can be estimated unbiasedly. The unbiased estimate was not monotone, however, and they showed how it can be improved by imposing the shape restriction.

Here we extend the isotonic analysis to estimate a regression function. The extension is motivated by the problem of estimating the velocity dispersions in astronomy which, in turn, is motivated by the problem of estimating the distribution of dark matter, as explained in [11]. Let $\mathbf{X} = (X_1, X_2, X_3)$ and $\mathbf{V} = (V_1, V_2, V_3)$ denote the three-dimensional position and velocity of a star and suppose that the distribution of **X** and **V** has a density that depends only on the Euclidean norms $\|\mathbf{x}\|$ and $\|\mathbf{v}\|$. We cannot observe a sample of three-dimensional star positions and velocities directly. What we can observe is the projected stellar positions and the line-of-sight of velocity components. With a proper choice of coordinates these become a sample of $(X_1, X_2, V_3)$. An important quantity to estimate is the velocity dispersion, $\sqrt{E(\|\mathbf{V}\|^2|R=r)}$, where $R = \sqrt{X_1^2 + X_2^2 + X_3^2}$. Due to the isotropy, $E(\|\mathbf{V}\|^2|R=r) = 3E(V_3^2|R=r)$. Clearly,

$$E(V_3^2|R=r) = \frac{\phi(r)}{f_\mathbf{r}(r)}, \tag{1}$$

where $\phi(r) = \int v_3^2 f_{\mathbf{r},v}(r,v)\,d\mathbf{v}$, $f_\mathbf{r}(x)$ is the density of $R$ and $f_{\mathbf{r},v}(r,v)$ is the joint density of $R$ and $\mathbf{V}$. We focus on estimating the function $\phi$ in (1). The function $f_\mathbf{r}$ can be estimated similarly (and was in [5]). Moreover, since positions are easier to measure than velocities, there may be additional information for estimating $f_\mathbf{r}$, as in [7].

Our main result is a version of the Kiefer–Wolfowitz theorem [8] comparing the empirical distribution function to its least concave majorant, but with a faster rate of convergence. In the next section, we relate $\phi$ to the derivative of a nondecreasing concave function $U$ that can be estimated unbiasedly by a nondecreasing, but nonconcave function $U_n^\#$. Letting $\tilde{U}_n$ be the least concave majorant of $U_n^\#$, we show that (under appropriate conditions) $\sup_t |\tilde{U}_n(t) - U_n^\#(t)| = O_p[\log(n)/n]$. In fact, a more general result is obtained, and the result is new even when specialized to the context of [5].

The paper is organized as follows. In Section 2 we establish the notation for several estimators and state the main result. The main result is useful in obtaining asymptotic distributions for estimators. Two such applications are presented in Section 3. An example with simulated data is presented in Section 4. Proofs are presented in Sections 5–8. Our study of Wicksell's problem follows Groeneboom and Jongbloed [5]. There are several other approaches to study Wicksell's problem in the literature. These include Hall and Smith's kernel method [6], Antoniadis, Fan and Gijbels's wavelet analysis [1] and Gobubev and Levit's asymptotically efficient estimation [4].

**2. The main results.** Let $\mathbf{X} = (X_1, X_2, X_3)$ be a random vector and $Z$ a nonnegative random variable. We suppose throughout that the joint distribution of $\mathbf{X}$ and $Z$ is invariant under orthogonal transformations of $\mathbf{X}$ and also that $\mathbf{X}$ and $Z$ have a joint density with respect to a product measure $\lambda^3 \times \nu$, where $\lambda^3$ is three-dimensional Lebesgue measure and $\nu$ is a sigma finite measure on the Borel sets of $[0,\infty)$. For example, if $\mathbf{X}$ is position and $Z = V_3^2$, then the invariance requires spatial symmetry. The relation between $E(\|\mathbf{V}\|^2|R=r)$ and $E(V_3^2|R=r)$ requires isotropy, but the estimator of $V_3^2$ does not. If $(\mathbf{X}, Z)$ are as described,

$$P[\mathbf{x} \leq \mathbf{X} \leq \mathbf{x} + d\mathbf{x}, z \leq Z \leq z + dz] = \rho(r; z)\, d\mathbf{x}\, \nu\{dz\}, \tag{2}$$

where $r^2 = x_1^2 + x_2^2 + x_3^2$ and $\rho$ is a nonnegative measurable function on $[0, \infty)^2$. Then $\rho(r) = \int_0^\infty \rho(r, z)\nu\{dz\}$ defines a marginal density for $\mathbf{X}$. Interest centers on estimating $E(Z|\mathbf{X} = \mathbf{x})$ from an observed sample of $(X_1, X_2, Z)$.

It is convenient to work with squared radii, $X = X_1^2 + X_2^2 + X_3^2$ and $Y = X_1^2 + X_2^2$, as in [5]. Let $f$ and $g$ denote the densities of the pairs $(X, Z)$ and $(Y, Z)$ with respect to $\lambda \times \nu$, where $\lambda$ is one-dimensional Lebesgue measure. Thus $f(x, z) = 2\pi \sqrt{x} \rho(\sqrt{x}, z)$ and

$$g(y, z) = \pi \int_y^\infty \frac{\rho(x, z)\, dx}{\sqrt{x - y}} = \frac{1}{2} \int_y^\infty \frac{f(x, z)\, dx}{\sqrt{x(x - y)}}.$$

Then

$$E(Z|X = x) = \frac{\varphi(x)}{\rho(\sqrt{x})},$$

where

$$\varphi(x) = \int z\rho(\sqrt{x}, z)\, dz. \tag{3}$$

We focus attention on estimating $\varphi$. Estimating $\rho(\sqrt{x})$ is then a special case with $Z \equiv 1$.

Let

$$\psi(y) = \int_0^\infty zg(y, z)\nu\{dz\}.$$

Then $\varphi$ and $\psi$ are related by

$$\psi(y) = \pi \int_y^\infty \frac{\varphi(x)\, dx}{\sqrt{x - y}}. \tag{4}$$

The transformation (4) can be inverted. Let

$$\Psi(y) = \int_y^\infty \frac{\psi(t)\, dt}{\sqrt{t - y}}. \tag{5}$$

Then

$$\Psi(y) = \pi^2 \int_y^\infty \varphi(x)\, dx \tag{6}$$

by reversing the orders of the integration in (5) and recognizing a Beta integral. It follows that $\varphi(x) = -\Psi'(x)/\pi^2$, so that estimation of $\varphi$ may proceed by estimating $\Psi$ and its derivative. Observe that, from equation (6) and the assumption that $Z$ is nonnegative, $\Psi$ is a nonincreasing function.

Suppose now that there is a sample $(X_{i,1}, X_{i,2}, Z_i)$, $i = 1, \ldots, n$, assumed to be i.i.d. as $(X_1, X_2, Z)$. Let $Y_i = X_{i,1}^2 + X_{i,2}^2$ and

$$\Psi_n^\#(y) = \frac{1}{n} \sum_{i:Y_i > y} \frac{Z_i}{\sqrt{Y_i - y}}. \tag{7}$$

Then $\Psi_n^\#(y)$ is an unbiased estimator of $\Psi(y)$ for each $y$, but is not monotone when viewed as a function of $y$: $\Psi_n^\#$ has an infinite jump at each observation $Y_i$, as indicated imperfectly in Figure 2. We call $\Psi_n^\#$ the *naive* estimator. This naive estimator can be improved by imposing shape restriction. If $\Psi_n^\#$ were square integrable, this could be accomplished by minimizing the integral of $(W - \Psi_n^\#)^2$ over non-increasing functions $W$, or equivalently, by minimizing

$$\int_0^\infty W^2(x)\,dx - 2\int_0^\infty W(x)\Psi_n^\#(x)\,dx. \tag{8}$$

The function $\Psi_n^\#$ is not square integrable, but it is integrable, so that (8) is well defined. Let $\tilde{\Psi}_n$ be the nonincreasing function $W$ that minimizes (8). Existence and uniqueness can be shown along the lines of Theorem 1.2.1 of Robertson, Wright and Dykstra [9], replacing the sums by integrals.

To describe $\tilde{\Psi}_n$ in more detail, let

$$U(x) = \int_0^x \Psi(t)\,dt = \int_0^\infty \int_0^\infty 2z[\sqrt{y} - \sqrt{(y-x)_+}]g(y,z)\nu\{dz\}\,dy, \tag{9}$$

where $z_+ = \max[0, z]$, and let

$$U_n^\#(x) = \int_0^x \Psi_n^\#(t)\,dt = \frac{1}{n} \sum_{i=1}^n 2Z_i(\sqrt{Y_i} - \sqrt{(Y_i - x)_+}). \tag{10}$$

Then $U$ is a nondecreasing, concave function, $U_n^\#$ a nondecreasing one and $\Psi_n^\#$ is the derivative of $U_n^\#$, at least *almost everywhere*. Let $\tilde{U}_n$ be the least concave majorant of $U_n^\#$. Then $\tilde{\Psi}_n$ is the right derivative of $\tilde{U}_n$. Letting

$$\varepsilon_n = \sqrt{\frac{\log(n)}{n}}, \tag{11}$$

the main results of the paper are:

THEOREM 2.1. *Suppose that $Z$ is a bounded random variable. Let $0 \leq t_0 < t_1 < \infty$ and suppose that $U$ is twice continuously differentiable on $[t_0, t_1]$ and that*

$$2\gamma_0 := \inf_{t_0 \leq t \leq t_1} [-U''(t)] > 0. \tag{12}$$

*Then*

$$\sup_{t_0 \leq t \leq t_1} |\tilde{U}_n(t) - U_n^\#(t)| = O_p[\varepsilon_n^2]. \tag{13}$$

THEOREM 2.2. *Suppose that $Z$ is a bounded random variable and $U$ is twice continuously differentiable near $x$ with $U''(x) < 0$. Then*

$$\sup_{|t-x| \leq \varepsilon_n} |\tilde{U}_n(t) - U_n^{\#}(t)| = o_p[\varepsilon_n^2]. \tag{14}$$

The proofs of Theorems 2.1 and 2.2 are presented in Sections 6 and 7.

## 3. Applications.

*Estimating $\Psi$.* From Theorem 2.1 $\tilde{U}_n$ and $U_n^{\#}$ are asymptotically equivalent. The following two theorems show that the derivatives $\Psi_n^{\#}$ and $\tilde{\Psi}_n$ are not. Let

$$\sigma^2(x) = \int_0^\infty z^2 g(x,z) \nu\{dz\}.$$

THEOREM 3.1. *For each $x \geq 0$ for which $\sigma^2(x)$ is continuous in a neighborhood of $x$,*

$$\sqrt{\frac{n}{\log n}}(\Psi_n^{\#}(x) - \Psi(x)) \Rightarrow N[0, \sigma^2(x)], \qquad \text{as } n \to \infty.$$

THEOREM 3.2. *For each $x > 0$ for which $\sigma^2(x)$ is continuous in a neighborhood of $x$,*

$$\sqrt{\frac{n}{\log n}}(\tilde{\Psi}_n(x) - \Psi(x)) \Rightarrow N\left[0, \frac{1}{2}\sigma^2(x)\right], \qquad \text{as } n \to \infty.$$

In the special case that $Z \equiv 1$, Theorems 3.1 and 3.2 are special cases of Theorems 2 and 5 in Groeneboom and Jongbloed [5]. We provide a different proof for Theorem 3.2 using Theorem 2.2 in Section 8.

*Estimating $\varphi$.* Extend $\Psi_n^{\#}$ and $\tilde{\Psi}_n$ to $(-\infty, \infty)$ by letting $\Psi_n^{\#}(t) = \Psi_n^{\#}(0)$ and $\tilde{\Psi}_n(t) = \tilde{\Psi}_n(0)$ for $t \leq 0$. Then, we may obtain a smooth estimator of $\Psi$ by using the kernel method,

$$\Psi_{n,s}(x) = \int_{-\infty}^{\infty} \frac{1}{b} K\left(\frac{t-x}{b}\right) \Psi_n^{\#}(t)\, dt, \tag{15}$$

where $K$ is a kernel and $b$ is a bandwidth. Due to the irregular behavior of $\Psi_n^{\#}$, $\Psi_{n,s}$ is not a nonincreasing function. If we replace the naive estimator by the isotonized estimator $\tilde{\Psi}_n$ in equation (15), we can obtain a smooth and nonincreasing estimator of $\Psi$,

$$\tilde{\Psi}_{n,s}(x) = \int_{-\infty}^{\infty} \frac{1}{b} K\left(\frac{t-x}{b}\right) \tilde{\Psi}_n(t)\, dt. \tag{16}$$

Differentiating (15) and (16) can give us the estimates of the derivative of the function $\Psi$. For instance, to estimate $\Psi'$, we have the estimators based on $\Psi^{\#}$ or $\tilde{\Psi}$,

$$\Psi'_{n,s}(x) = -\int_{-\infty}^{\infty} \frac{1}{b^2} K'\left(\frac{t-x}{b}\right) \Psi_n^{\#}(t)\, dt \tag{17}$$

and

$$\tilde{\Psi}'_{n,s}(x) = -\int_{-\infty}^{\infty} \frac{1}{b^2} K'\left(\frac{t-x}{b}\right) \tilde{\Psi}_n(t)\, dt. \tag{18}$$

Hall and Smith [6] studied $\Psi'_{n,s}$ for the density estimation problem ($Z \equiv 1$). They showed that the optimal order for the bandwidth is $n^{-1/6}$ under broad conditions that also imply that $n^{1/6}[\Psi'_{n,s}(t) - \Psi'(t)]$ is asymptotically normal. As a consequence of Theorem 2.1, we may show that $\Psi'_{n,s}$ and $\tilde{\Psi}'_{n,s}$ have the same asymptotic distribution, under broad conditions that include the optimal order for the bandwidth.

THEOREM 3.3. *Suppose that $K$ is twice continuously differentiable and has support $[-1, 1]$. If $Z$ is bounded, $\varphi$ is positive and continuous at $0 < x < \infty$, $b = b_n \to 0$ and $nb \to \infty$, then*

$$|\tilde{\Psi}'_{n,s}(x) - \Psi'_{n,s}(x)| = O_p\left[\frac{\log(n)}{nb^2}\right].$$

PROOF. Since $\varphi$ is positive and continuous at $x$, there are $t_0 < x < t_1$ for which (12) holds. If $n$ is sufficiently large, then $[x - b, x + b] \subseteq [t_0, t_1]$ and, therefore,

$$\begin{aligned}
|\tilde{\Psi}'_{n,s}(x) - \Psi'_{n,s}(x)| &= \left|\int_{x-b}^{x+b} \frac{1}{b^2} K'\left(\frac{t-x}{b}\right)[\tilde{\Psi}_n(t) - \Psi_n^{\#}(t)]\, dt\right| \\
&\leq \frac{1}{b^3} \int_{x-b}^{x+b} \left|K''\left(\frac{t-x}{b}\right)\right| \left|\tilde{U}_n(t) - U_n^{\#}(t)\right| dt \\
&= O_p\left[\frac{\log(n)}{nb^2}\right]. \qquad \square
\end{aligned}$$

**4. Simulations.** To illustrate the nature of the estimators, we use a simulated sample from Plummer's distribution [2], page 205. In Plummer's distribution, the joint density of $\mathbf{X} = (X_1, X_2, X_3)$ and $Z = V_3^2$ is

$$\frac{c_0}{\beta^5 \sqrt{z}}\left[\frac{\beta}{\sqrt{1 + (1/3)r^2}} - \frac{1}{2} z\right]_+^5,$$

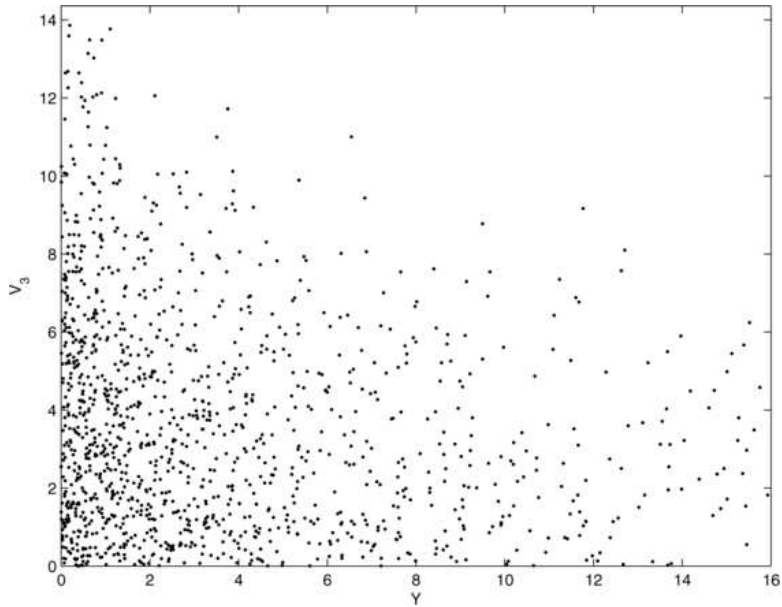

Fig. 1. *Simulated $Y$ and $V_3$ pair.*

where $c_0$ is a normalizing constant and $\beta$ is a parameter that is related to the velocity dispersion through $E[Z|R=r] = \beta/[6\sqrt{1+(1/3)r^2}]$. For this case,

$$\Psi(y) = \frac{\sqrt{3}\pi\beta}{48} \frac{1}{(1+(1/3)y)^2}.$$

A sample of 1500 $(Y, V_3)$ pairs, simulated from the above distributions with $\beta = 200$, is shown in Figure 1. The naive estimator $\Psi_n^\#$ may be computed from these data using equation (7). This function is shown in Figure 2 along with the improved estimator $\tilde{\Psi}_n$ and kernel estimator $\tilde{\Psi}_{n,s}$. As expected, that $\Psi_n^\#$ is a highly irregular function and $\tilde{\Psi}_n$ is a decreasing step function. The estimator $\tilde{\Psi}_{n,s}$ is a smooth function. Figure 3 shows the estimator $\tilde{\Psi}'_{n,s}$ of the first derivative of $\Psi$. In all cases, the dashed line represents the true function computed directly from the true distribution. The bandwidths in Figures 2 and 3 were 1.5 and 3.7, respectively. These were chosen by inspection to compromise between fit and smoothness.

**5. Localization lemmas.** A question that arises in the proofs of Theorems 2.1 and 2.2 is the relation between the least concave majorant of a restricted function and the restriction of the least concave majorant. This question is of some independent technical interest.

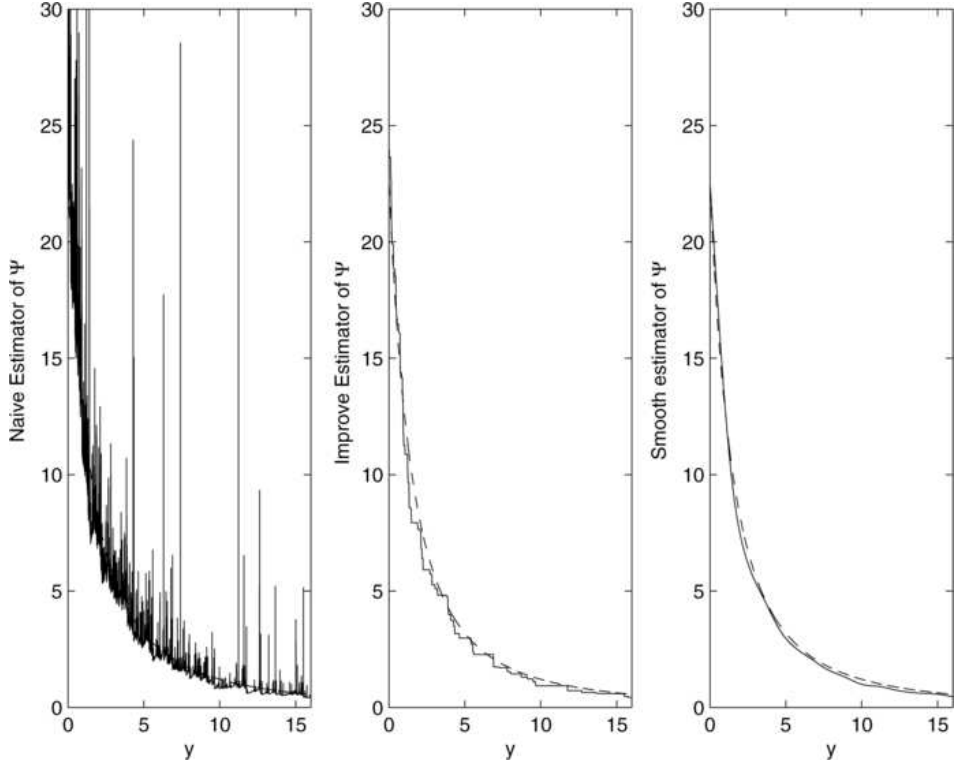

FIG. 2. *Naive, isotonized and smooth estimates of $\Psi$.*

LEMMA 5.1. *Let $f$ be a bounded continuous function on $[0,\infty)$, let $\tilde{f}$ be its least concave majorant, and let $0 < a < a' < \infty$ and $c = (a+a')/2$. If $f(c) > [\tilde{f}(a) + \tilde{f}(a')]/2$, then $f(t) = \tilde{f}(t)$ for some $a \leq t \leq a'$.*

PROOF. If $f(c) > [\tilde{f}(a) + \tilde{f}(a')]/2$, then clearly $\tilde{f}(c) > [\tilde{f}(a) + \tilde{f}(a')]/2$. Let $g(t) = \tilde{f}(t)$ for $t \leq a$ and $t \geq a'$, and let

$$g(t) = (1-\varepsilon)\tilde{f}(t) + \varepsilon\left[\frac{(t-a)\tilde{f}(a') + (a'-t)\tilde{f}(a)}{a'-a}\right]$$

for $a \leq t \leq a'$, where $0 < \varepsilon < 1$. Then $g$ is concave on $[0,\infty)$ for any $0 < \varepsilon < 1$, and $g(t) \leq \tilde{f}(t)$ for all $t$ with strict inequality if $t = c$. If $f(t) < \tilde{f}(t)$ for $a \leq t \leq a'$, then there is an $\varepsilon > 0$ for which $f(t) < g(t)$ for $a \leq t \leq a'$, and this would contradict the definition of $\tilde{f}$. □

LEMMA 5.2. *Let $f$ be a bounded function on $[0,\infty)$ and let $\tilde{f}$ be its least concave majorant. Further, let $0 \leq z_0 \leq x_0 \leq t_0 < t_1 < x_1 < z_1 < \infty$ with strict inequality throughout if $t_0 > 0$, and let $f^*$ be the least concave*

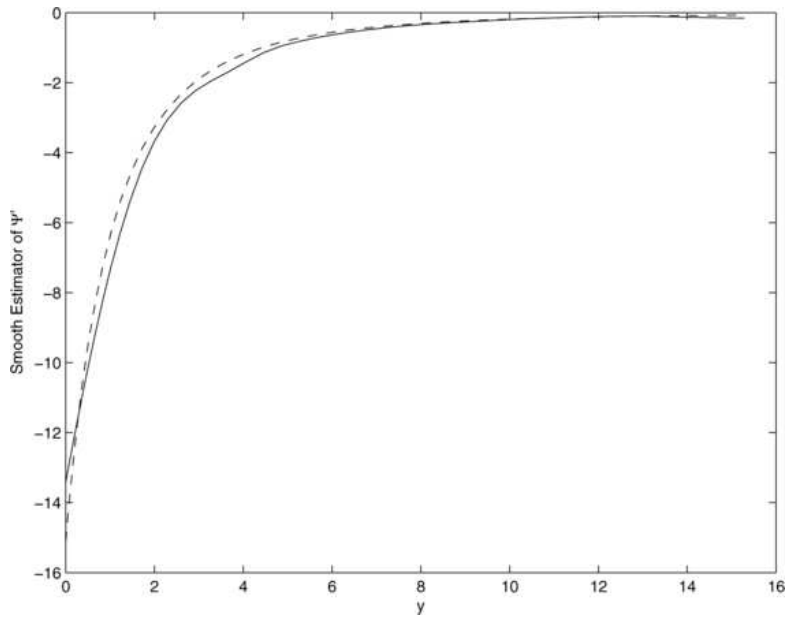

Fig. 3. *Estimator of $\Psi'$*.

*majorant of the restriction of $f$ to $[z_0, z_1]$. If $\tilde{f}(x_i) = f(x_i)$, $i = 0, 1$, then $\tilde{f}(t) = f^*(t)$ for $x_0 \leq t \leq x_1$.*

PROOF. It is clear that $f^*(t) \leq \tilde{f}(t)$ for all $t \in [x_0, x_1]$ and suffices to establish the reverse inequality. It is also clear that $\tilde{f}(x_i) = f^*(x_i)$, since $f^*(x_i) \geq f(x_i)$. Let

$$f^+(t) = \begin{cases} \tilde{f}(t), & \text{if } t \leq x_0, \\ f^*(t), & \text{if } x_0 \leq t \leq x_1, \\ \tilde{f}(t), & \text{if } t \geq x_1. \end{cases}$$

Then $f^+ \geq f$ everywhere. It is shown below that $f^+$ is concave. It then follows that $f^+ \geq \tilde{f}$, and this implies $f^*(t) \geq \tilde{f}(t)$ for $x_0 \leq t \leq x_1$.

The restrictions of $f^+$ to each of the intervals $[0, x_0]$, $[x_0, x_1]$ and $[x_1, \infty)$ are clearly concave. So, the question is whether the derivative of $f^+$ is nonincreasing as $t$ crosses $x_0$ and $x_1$. Consider $x_0$. The issue does not arise if $t_0 = 0$; and if $t_0 > 0$, then the left- and right-hand derivatives of $f^+$ at $x_0$ are related by

$$D_r f^+(x_0) = \lim_{\varepsilon \downarrow 0} \frac{f^+(x_0 + \varepsilon) - f^+(x_0)}{\varepsilon}$$
$$\leq \lim_{\varepsilon \downarrow 0} \frac{\tilde{f}(x_0 + \varepsilon) - \tilde{f}(x_0)}{\varepsilon}$$

$$= D_r \tilde{f}(x_0) \leq D_\ell \tilde{f}(x_0) = D_l f^+(x_0).$$

Global concavity follows from this and a related calculation at $x_1$. □

**6. Proof of Theorem 2.1.** The assumptions of Theorem 2.1 are made throughout this section, and the proof is divided into several steps. Let $0 \leq t_0 < t_1 < \infty$ and $\gamma_0$ be as in the statement of that theorem.

STEP 1. For suitably chosen $w_0 \leq t_0$, with equality when $t_0 = 0$, and $w_1 > t_1$, let $U_n^*$ be the least concave majorant of the restriction of $U_n^\#$ to $[w_0, w_1]$. Then

(19) $$\lim_{n \to \infty} P[\tilde{U}_n(t) = U_n^*(t) \text{ for all } t \in [t_0, t_1]] = 1.$$

STEP 2. Next, partition the interval $[w_0, w_1]$ into $k_n$ equally-spaced subintervals $[a_j^{k_n}, a_{j+1}^{k_n}]$, $j = 0, 1, \ldots, k_n - 1$, where $k_n$ is a sequence of positive integers for which $\log(n) \leq k_n \leq \sqrt{n}$, to be specified in detail later. Thus,

$$a_j^{k_n} = w_0 + \frac{wj}{k_n}, \qquad j = 0, 1, \ldots, k_n,$$

where $w = w_1 - w_0$. Define two continuous piecewise linear functions $L^{k_n}$ and $L_n^{k_n}$ with knots at $a_j^{k_n}$ by $L^{k_n}(a_j^{k_n}) = U(a_j^{k_n})$ and $L_n^{k_n}(a_j^{k_n}) = U_n^\#(a_j^{k_n})$ for $j = 0, 1, \ldots, k_n$. Let $A_n$ be the event that $L_n^{k_n}$ is concave. In Step 2, it is shown that $P(A_n) \to 1$ as $n \to \infty$, for suitably chosen $k_n$.

Let $\|f\| = \sup\{|f(t)| : t \in [w_0, w_1]\}$.

STEP 3. If $A_n$ occurs, then $\|U_n^* - U_n^\#\| \leq 2\|L_n^{k_n} - U_n^\# + U - L^{k_n}\| + 2\|L^{k_n} - U\|$.

STEP 4. There are constants $C$, $\Gamma_1$ and $\Gamma_2$ for which

(20) $$\|L^{k_n} - U\| \leq \frac{\Gamma_1}{k_n^2}$$

and

(21) $$P\left[\|L_n^{k_n} - U_n^\# - (U - L^{k_n})\| > \frac{\sqrt{\lambda \log(n)}}{\sqrt{n} k_n}\right]$$
$$< C k_n^3 [\log(n)]^7 \exp\left[-\frac{\lambda \log(n)}{\Gamma_2}\right]$$

for all $0 < \lambda \leq \log(n)$ for all sufficiently large $n$.

It is then shown that the results in Steps 1–4 imply (13) with suitable choice of $k_n$.

In the first proposition, let $w_0 \leq t_0$ and $w_1 > t_1$ be so chosen that $\sup_{w_0 \leq t \leq w_1} U''(t) \leq -\gamma_0$ and $w_0 < t_0$ if $t_0 > 0$. In the proof, use is made of the following elementary observations: $\sup_{0 \leq t < \infty} |U_n^\#(t) - U(t)| \to 0$ w.p.1 by the strong law of large numbers applied for fixed $t$ and the monotonicity of $U_n^\#$, as in Polya's theorem; and then

$$(22) \quad \sup_{0 \leq t < \infty} |\tilde{U}_n(t) - U(t)| \leq \sup_{0 \leq t < \infty} |U_n^\#(t) - U(t)| \to 0 \quad \text{w.p.1}$$

by Marshall's lemma.

PROPOSITION 6.1. *Relation* (19) *holds.*

PROOF. By Lemma 5.2, it suffices to show that $\tilde{U}_n(t) = U_n^\#(t)$ for some $t \in [w_0, t_0]$ and also for some $t \in [t_1, w_1]$, both with probability approaching one. The two cases are similar. Details are supplied for the first. The relation is clear if $t_0 = 0$. If $t_0 > 0$, let $x_0 = (w_0 + t_0)/2$. Then it suffices to show that $U_n^\#(x_0) > [\tilde{U}_n(w_0) + \tilde{U}_n(t_0)]/2$ with probability approaching one, by Lemma 5.1. This is clear, however, from (22) and the strict concavity of $U$ which implies $U(x_0) > [U(w_0) + U(t_0)]/2$. □

This accomplishes Step 1. For the next step, let $Q$ denote the distribution of $(Y, Z)$, $\|\cdot\|_r$ the norm in $L^r(Q)$ and $Q_n^\#$ the empirical distribution of $(Y_i, Z_i)$, $i = 1, \ldots, n$. Then

$$U_n^\#(t) = Q_n^\# h_t := \frac{1}{n} \sum_{i=1}^n h_t(Y_i, Z_i),$$

where

$$h_t(y, z) = 2z[\sqrt{y} - \sqrt{(y-t)_+}]$$

and $(y - t)_+ = \max[y - t, 0]$. For later reference observe that

$$(23) \quad |h_t - h_s| \leq 2c_0 \sqrt{|t - s|},$$

where $c_0$ is an upper bound for $Z$.

LEMMA 6.1. *Let $\omega_0 \leq a < b \leq \omega_1$. There is a constant $\Gamma_3$ such that if $|s - t| < 1$, then*

$$(24) \quad \|h_t - h_s\|_2 \leq \Gamma_3 |t - s| \sqrt{\log \frac{2}{|t - s|}},$$

*and for all $a \leq t \leq b$,*

$$(25) \quad \left\| h_t - \frac{(t-a)h_b + (b-t)h_a}{b-a} \right\|_2 \leq \Gamma_3 (b - a).$$

PROOF. The square of the left-hand side of (24) is at most

$$4c_0^2 \int_s^{t+1} [\sqrt{(y-s)} - \sqrt{(y-t)_+}]^2 g(y)\, dy$$
$$+ 4c_0^2 \int_{t+1}^\infty [\sqrt{(y-s)} - \sqrt{(y-t)}]^2 g(y)\, dy.$$

Letting $c_1$ be an upper bound for $g$, the first integral is at most $4c_0^2 c_1 (s-t)^2 \int_0^{1/|s-t|} (\sqrt{z} - \sqrt{(z-1)_+})^2\, dz \leq 4c_0^2 c_1 (s-t)^2 \log(1/|s-t|)$ and the second integral is less than $4c_0^2 c_1 (s-t)^2$. Thus (24) follows. The square of the left-hand side of (25) is at most

$$4c_0^2 \int_0^\infty \left[\sqrt{(y-t)_+} - \frac{(t-a)\sqrt{(y-b)_+} + (b-t)\sqrt{(y-a)_+}}{b-a}\right]^2 g(y)\, dy$$

for $a \leq t \leq b$. Letting $y = a + (b-a)z$ and $t = a + (b-a)s$, the last integral is at most

$$4c_0^2 c_1 (b-a)^2 \int_0^\infty \{\sqrt{(z-s)_+} - [s\sqrt{(z-1)_+} + (1-s)\sqrt{z}]\}^2\, dz.$$

The integrand in the last expression is continuous and of order $1/z^3$ as $z \to \infty$ uniformly in $0 \leq s \leq 1$. The lemma follows. □

PROPOSITION 6.2. *There is a constant $\gamma_1 > 0$, for which*

$$P[L_n^{k_n} \text{ is concave }] \geq 1 - k_n \exp\left[-\frac{\gamma_1 n}{k_n^2}\right].$$

PROOF. Since $L_n^{k_n}$ is piecewise linear and $a_j^{k_n}$ are equally-spaced, the event that $L_n^{k_n}$ is concave is

$$A_n = \bigcap_{j=2}^{k_n-1} \{U_n^\#(a_{j+1}^{k_n}) + U_n^\#(a_{j-1}^{k_n}) - 2U_n^\#(a_j^{k_n}) \leq 0\}.$$

For a fixed $j$, let $h_{n,j} = h_{a_{j+1}^{k_n}} + h_{a_{j-1}^{k_n}} - 2h_{a_j^{k_n}}$. Then $U_n^\#(a_{j+1}^{k_n}) + U_n^\#(a_{j-1}^{k_n}) - 2U_n^\#(a_j^{k_n}) = Q_n^\# h_{n,j}$, and

$$P[U_n^\#(a_{j+1}^{k_n}) + U_n^\#(a_{j-1}^{k_n}) - 2U_n^\#(a_j^{k_n}) \geq 0] = P[Q_n^\# h_{n,j} - Qh_{n,j} \geq -Qh_{n,j}].$$

Here $Qh_{n,j} = U(a_{j+1}^{k_n}) + U(a_{j-1}^{k_n}) - 2U(a_j^{k_n})$ and $Qh_{n,j} \leq -\gamma_0 w^2/k_n^2$. Further, $|h_{n,j}| \leq 2c_0\sqrt{w/k_n}$ by (23), and the variance of $h_{n,j}(Y,Z)$ is at most $2\Gamma_3^2 w^2/k_n^2$ by the previous lemma. So, the last displayed expression is at most

$$P\left[Q_n^\# h_{n,j} - Qh_{n,j} \geq \frac{\gamma_0 w^2}{k_n^2}\right] \leq \exp\left[-\frac{\gamma_0^2 w^2 n}{4\Gamma_3 k_n^2 + 2c_0 w^{1/2}\gamma_0 k_n^2}\right]$$

by Bernstein's inequality, [10], page 102. The lemma follows with $\gamma_1 = \gamma_0^2 w^2 / (4\Gamma_3^2 + 2c_0\gamma_0\sqrt{w})$ by summing over $j$. □

The previous two lemmas accomplish Step 2. The next step is similar to Lemma 5 in [8]. Recall Marshall's lemma: If $h$ is a concave function, then $\|U_n^* - h\| \le \|U_n^\# - h\|$.

PROPOSITION 6.3. *If $A_n$ occurs, then $\|U_n^* - U_n^\#\| \le 2\|L_n^{k_n} - U_n^\# + U - L^{k_n}\| + 2\|L^{k_n} - U\|$.*

PROOF. If $A_n$ occurs, then $L_n^{k_n}$ is a concave function and, therefore, $U_n^* - U_n \le \|U_n^* - L_n^{k_n}\| + \|L_n^{k_n} - U_n\| \le 2\|L_n^{k_n} - U_n\| \le 2\|L_n^{k_n} - U_n + U - L^{k_n}\| + 2\|L^{k_n} - U\|$. □

For the final step, let

$$G_n f = \sqrt{n}(Q_n^\# f - Qf) = \frac{1}{\sqrt{n}} \sum_{i=1}^{n} [f(Y_i, Z_i) - Qf]$$

for $f \in L^1(Q)$, and $\|G_n\|_\mathcal{F} = \sup_{f \in \mathcal{F}} |G_n f|$ for $\mathcal{F} \subseteq L^1(Q)$.

PROPOSITION 6.4. *There are constants $C$, $\Gamma_1 > 0$, and $\Gamma_2$, for which (20) and (21) hold for all sufficiently large $n$.*

PROOF. That (20) holds follows from a simple Taylor series expansion, as in [8].

The inequality (21) requires more effort. First observe that

$$\|L_n^{k_n} - U_n^\# - (L^{k_n} - U)\| = \max_{1 \le j \le k_n} \max_{a_{j-1}^{k_n} \le t \le a_j^{k_n}} |L_n^{k_n}(t) - U_n^\#(t) - [L^{k_n}(t) - U(t)]|.$$

For fixed $n$ and $j$, let

$$h_t^* = h_t - \frac{(t - a_{j-1}^{k_n})h_{a_j^{k_n}} + (a_j^{k_n} - t)h_{a_{j-1}^{k_n}}}{a_j^{k_n} - a_{j-1}^{k_n}}$$

and

$$f_s = \frac{\sqrt{k_n}}{2c_0} h_{a_{j-1}^{k_n} + s/k_n}^*$$

for $a_{j-1}^{k_n} \le t \le a_j^{k_n}$ and $0 \le s \le 1$. Then $|f_s| \le 1$, and the variance of $f_s(Y, Z)$ is at most $\sigma^2 := \Gamma_3^2 w^2 / (4c_0^2 k_n)$ by (23) and Lemma 6.1. Moreover,

$$L_n^{k_n}(t) - U_n^\#(t) - [L^{k_n}(t) - U(t)] = Q_n^\# h_t^* - Q h_t^*,$$

$$\max_{a_{j-1}^{k_n} \le t \le a_j^{k_n}} |Q_n^\# h_t^* - Q h_t^*| = \max_{0 \le t \le 1} \frac{2c_0}{\sqrt{nk_n}} G_n f_t,$$

and there is a constant $C$ for which

$$P[\|G_n\|_{\mathcal{F}} > \tau] \leq \frac{C}{\sigma^4}\left[1 \vee \frac{\tau}{\sigma}\right]^7 \exp\left[-\frac{1}{2}\frac{\tau^2}{\sigma^2 + (3+\tau)/\sqrt{n}}\right]$$

for all $\tau > 0$ by Theorem 2.14.16 of [10]. The condition (2.14.6) of [10] is satisfied with $V = 2$ here by (23). The constant $C$ here does not depend on $n$ or $j$. Thus, combining the last three displayed expressions, letting

$$\tau = \frac{1}{2c_0}\sqrt{\frac{\lambda \log(n)}{k_n}}$$

and summing over $j$ leads to (21), with $\Gamma_2 = 2(\Gamma_3^2 w^2 + 12c_0)$. $\square$

PROOF OF THEOREM 2.1. Let $k_n$ be the least integer that exceeds $[\gamma_1 n / 4\log n]^{1/2}$. Then $k_n = O(n/\log n)^{1/2}$ and $P(A_n) \geq 1 - n^{-2}$ for all sufficiently large $n$ by Proposition 6.2. Choose $\lambda > 4\Gamma_2$ in (21). Let $D_n$ be the event $\tilde{U}_n(t) = U_n^*(t)$ for all $t \in [t_0, t_1]$. Then $A_n \cap D_n$ implies

$$\sup_{t_0 \leq t \leq t_1} |\tilde{U}_n - U_n^\#| = \sup_{t_0 \leq t \leq t_1} |U_n^* - U_n^\#| \leq \frac{\Gamma_1}{k_n^2} + \frac{\sqrt{4\lambda}\log n}{\sqrt{\gamma_1}n} = O(n^{-1}\log n),$$

(26)

by Propositions 6.1, 6.3 and 6.4. So, for large constant $M$,

$$P\left(\sup_{t_0 \leq t \leq t_1} |\tilde{U}_n - U_n^\#| > M\varepsilon_n^2\right) \leq P(A_n') + P(D_n'),$$

establishing (13). $\square$

**7. Proof of Theorem 2.2.** Theorem 2.2 can be proved by modifying the proof of Theorem 2.1. The following simpler and more transparent proof was suggested by a referee.

PROOF OF THEOREM 2.2. Suppose that $t_0 > 0$ and for $|s| \leq 1$, let

$$V_n(s) = \frac{(s+1)U_n^\#(t_0 + \varepsilon_n) + (1-s)U_n^\#(t_0 - \varepsilon_n)}{2}$$
$$+ U(t_0 + \varepsilon_n s) - \frac{(s+1)U(t_0 + \varepsilon_n) + (1-s)U(t_0 - \varepsilon_n)}{2}$$

and

$$R_n(s) = \frac{1}{\varepsilon_n^2}[U_n^\#(t_0 + \varepsilon_n s) - V_n(s)].$$

Observe that $V_n$ is concave and $R_n$ continuous, both on $[-1, 1]$. After some algebra, the process $R_n$ may be rewritten as

$$R_n(s) = \frac{1}{n\varepsilon_n^2}\sum_{i=1}^n [f_s(Y_i, Z_i) - Ef_s(Y_i, Z_i)],$$

where
$$f_s = h_{t_0+\varepsilon_n s} - \frac{(s+1)h_{t_0+\varepsilon_n} + (1-s)h_{t_0-\varepsilon_n}}{2}.$$

By Lemma 6.1,
$$\|f_t\|_2 \leq 2\Gamma_3 \varepsilon_n \quad \text{and} \quad \|f_t - f_s\|_2 \leq \Gamma_3 \varepsilon_n \sqrt{\log\left(\frac{2}{\varepsilon_n}\right)} |t-s|$$

for $-1 \leq s, t \leq 1$. It follows that
$$E[R_n(t)^2] \leq \frac{2\Gamma_3^2}{n\varepsilon_n^4}\varepsilon_n^2 \leq \frac{2\Gamma_3^2}{\log(n)}$$

and
$$E[R_n(t) - R_n(s)]^2 \leq \frac{\Gamma_3^2}{n\varepsilon_n^4}\varepsilon_n^2 \log\left(\frac{2}{\varepsilon_n}\right)(t-s)^2 \log\left(\frac{1}{|t-s|}\right)$$
$$\leq \Gamma_3^2 (t-s)^2 \log\left(\frac{1}{|t-s|}\right).$$

So, $R_n$ are tight in $C[-1,1]$, $R_n(s) \to^p 0$ for each $|s| \leq 1$, and therefore, $\sup_{|s|\leq 1}|R_n(s)| \to^p 0$. It follows that
$$\sup_{|s|\leq 1} |U_n^\#(t_0 + \varepsilon s) - V_n(s)| = o_p(\varepsilon_n^2).$$

So, $\sup_{|s|\leq 1}|U_n^*(t_0 + \varepsilon s) - V_n(s)| \leq \sup_{|s|\leq 1}|U_n^\#(t_0 + \varepsilon s) - V_n(s)|$, and the theorem follows from $|U_n^\# - \tilde{U}_n| \leq |U_n^\# - V_n| + |\tilde{U}_n - V_n|$.

## 8. Proof of Theorems 3.1 and 3.2.

PROOF OF THEOREM 3.1. For each $x \geq 0$, let
$$\Psi_i = \frac{Z_i}{\sqrt{Y_i - x}} 1_{\{Y_i > x\}}, \qquad i = 1, \ldots, n.$$

Then, for each $c > 0$,
$$P[\Psi_i > c] = \int_0^\infty \int_x^{x+z^2/c^2} g(y,z) \, dy \, \nu\{dz\}.$$

So
$$(27) \qquad c^2 P[\Psi_i > c] = \int_0^1 \int_0^\infty z^2 g\left(x + \frac{t^2 z^2}{c^2}, z\right) \nu\{dz\} \, dt \to \sigma^2(x)$$

as $c \to \infty$, if $\sigma^2$ is continuous in a neighborhood of $x$, using Pratt's theorem to justify the interchange of limit and integral. Theorem 3.1 follows easily from [3], Theorem 4, page 323. □

PROOF OF THEOREM 3.2. To begin fix $0 < x < \infty$. From Theorem 2.2 and more basic considerations, there is a sequence $0 < \delta_n \leq \varepsilon_n$ for which $\delta_n = o(\varepsilon_n)$, $\log(\delta_n^{-1}) \sim \log(\varepsilon_n^{-1})$ and

$$\sup_{|x'-x|\leq \varepsilon_n} |\tilde{U}_n(x') - U_n^{\#}(x')| = o_p(\delta_n^2). \tag{28}$$

Since $\tilde{U}_n$ is concave,

$$\frac{\tilde{U}_n(x+\delta_n) - \tilde{U}_n(x)}{\delta_n} \leq \tilde{\Psi}_n(x) \leq \frac{\tilde{U}_n(x) - \tilde{U}_n(x-\delta_n)}{\delta_n}. \tag{29}$$

From (28),

$$\frac{\tilde{\Psi}_n(x) - \Psi(x)}{\varepsilon_n} \geq \frac{U_n^{\#}(x+\delta_n) - U_n^{\#}(x) - \delta_n \Psi(x)}{\delta_n \varepsilon_n} - o_p(1),$$

$$\frac{\tilde{\Psi}_n(x) - \Psi(x)}{\varepsilon_n} \leq \frac{U_n^{\#}(x-\delta_n) - U_n^{\#}(x) - \delta_n \Psi(x)}{\delta_n \varepsilon_n} + o_p(1).$$

Here

$$E\left[\frac{U_n^{\#}(x \pm \delta_n) - U_n^{\#}(x)}{\delta_n}\right] = \Psi(x) + O(\delta_n).$$

Next, consider the asymptotic variance of $U_n^{\#}(x+\delta_n) - U_n^{\#}(x)$,

$$\mathrm{Var}[U_n^{\#}(x+\delta_n) - U_n^{\#}(x)] = \frac{1}{n}\mathrm{Var}[2Z_i(\sqrt{[Y_i - x]_+} - \sqrt{[Y_i - x - \delta_n]_+})]$$

$$= \frac{1}{n}E[4Z_i^2(Y_i - x)1_{\{x < Y_i < Y_i + \delta_n\}}]$$

$$+ \frac{1}{n}E[4Z_i^2(\sqrt{Y_i - x} - \sqrt{Y_i - x - \delta_n})1_{\{Y_i > x + \delta_n\}}]$$

$$+ \frac{1}{n}(U(x+\delta_n) - U(x))^2.$$

It is easy to see that

$$E[4Z_i^2(Y_i - x)1_{\{x < Y_i < x + \delta_n\}}] + \frac{1}{n}(U(x+\varepsilon_n) - U(x))^2 = O(\delta_n^2),$$

provided that $\int z^2 g(x, z)\, dz$ is continuous at a neighborhood of $x$, and

$$E[4Z_i^2(\sqrt{Y_i - x} - \sqrt{Y_i - x - \delta_n})1_{\{Y_i > x + \delta_n\}}] \sim -\delta_n^2 \log \delta_n \int z^2 g_{\mathbf{y},z}(x, z)\, dz$$

under the same conditions. So

$$\mathrm{Var}\left[\frac{U_n^{\#}(x+\delta_n) - U_n^{\#}(x)}{\delta_n}\right] = \frac{1}{n}\log\frac{1}{\delta_n}\sigma^2(x) \sim \frac{1}{2}\varepsilon_n^2 \sigma^2(x).$$

A dual result is easily obtained from $U_n^{\#}(x) - U_n^{\#}(x - \delta)$. So

$$\frac{U_n^{\#}(x \pm \delta_n) - U_n^{\#}(x) - \delta_n \Psi(x)}{\delta_n \varepsilon_n} \Rightarrow N\left(0, \frac{1}{2}\sigma^2(x)\right)$$

from the Lindeberg central limit theorem for triangular arrays. $\square$

**Acknowledgments.** We would like to thank the Editor, an Associate Editor and two referees for their constructive comments and suggestions which helped improve our presentation greatly.

DEPARTMENT OF MATHEMATICS AND STATISTICS
UNIVERSITY OF MARYLAND BALTIMORE COUNTY
BALTIMORE, MARYLAND 21250
USA
E-MAIL: wangxiao@umbc.edu

DEPARTMENT OF STATISTICS
UNIVERSITY OF MICHIGAN
ANN ARBOR, MICHIGAN 48109–1092
USA
E-MAIL: michaelw@umich.edu
URL: www.stat.lsa.umich.edu/~michaelw